\documentclass[11pt]{article}
\usepackage{graphicx}
\usepackage{amssymb}
\usepackage{amsmath}
\usepackage{amsthm}
\usepackage{subcaption}
\usepackage{tikz}
\usetikzlibrary{arrows, decorations.markings, decorations.pathmorphing, backgrounds, positioning, fit, petri}
\usetikzlibrary{shapes.geometric}

\newtheorem{thm}{Theorem}[section]
\newtheorem{prop}[thm]{Proposition} 
\newtheorem{cor}[thm]{Corollary} 
\newtheorem{lemma}[thm]{Lemma}

\theoremstyle{definition}
\newtheorem{example}{Example}
\newtheorem{defn}[thm]{Definition}
\newtheorem{remark}{Remark}

\begin{document}

\title{Relative Non-Positive Immersion}
\author{Jens Harlander, Stephan Rosebrock}
\maketitle

\begin{abstract}
A 2-complex $K$ has collapsing non-positive immersion if for every combinatorial immersion $X\to K$, where $X$ is finite, connected and does not allow collapses, either $\chi(X)\le 0$ or $X$ is point. This concept is due to Wise who also showed that this property implies local indicability of the fundamental group $\pi_1(K)$. In this paper we study a relative version of collapsing non-positive immersion that can be applied to 2-complex pairs $(L,K)$: The pair has relative collapsing non-positive immersion if for every combinatorial immersion $f\colon X\to L$, where $X$ is finite, connected and does not allow collapses, either $\chi(X)\le \chi(Y)$, where $Y$ is the essential part of the preimage $f^{-1}(K)$, or $X$ is a point. We show that under certain conditions a transitivity law holds: If $(L,K)$ has relative collapsing non-positive immersion and $K$ has collapsing non-positive immersion, then $L$ has collapsing non-positive immersion. This article is partly motivated by the following open question: Do reduced injective labeled oriented trees have collapsing non-positive immersion? We answer this question in the affirmative for certain important special cases.
\end{abstract}



\section{Introduction}
Suppose $\mathcal P$ is a property that can be attributed to a 2-complex. A relative version of $\mathcal P$ which can be applied to 2-complex pairs $(L,K)$, where $K$ is a subcomplex of $L$, can be a powerful tool, especially if a transitivity law holds: If $(L,K)$ has property $\mathcal P$ relative to $K$ and $K$ has property $\mathcal P$, then $L$ has property $\mathcal P$. A good example of a property that has a relative version and satisfies a transitivity law is \lq\lq vanishing reduced homology\rq\rq. In \cite{HarRose2017} and \cite{HarRose2020} we used relative versions of combinatorial asphericity to prove that injective labeled oriented trees are aspherical. This answered the long standing asphericity question for ribbon disc complements in the alternating case. In this paper we study relative non-positive immersion. 

An $n$-cell in a complex is called {\em free} if it is used exactly once in an attaching map of a unique $(n+1)$-cell. A cell that contains a free cell in its boundary can be collapsed without changing the homotopy type of the complex. A 2-complex map $f\colon X\to K$ is {\em combinatorial} if it maps open cells homeomorphically to open cells. It is an {\em immersion} if it is locally injective.

\begin{defn} A 2-complex $K$ has {\em collapsing non-positive immersion} if for every combinatorial immersion $X\to K$, where $X$ is finite, connected, and without free vertices or edges, either $\chi(X)\le 0$ or $X$ is a point.
\end{defn}

A number of variants of non-positive immersion were defined by Wise. See Definition 1.2 in \cite{Wise22}. See also \cite{Wise04}. Among many other things Wise showed that if $K$ has non-positive immersion then $\pi_1(K)$ is locally indicable in case that this group is not trivial.

The {\em essential part} of a 2-complex consists of all the 2-cells together with the edges that are used in the attaching paths for the 2-cells. 

\begin{defn}
Let $L$ be a 2-complex and $K$ be a subcomplex. Then $(L,K)$ has {\em relative collapsing non-positive immersion} if for every combinatorial immersion $f\colon X\to L$, where $X$ is finite, connected, and without free vertices or edges, either $\chi(X)\le \chi(Y)$, where $Y$ is the essential part of the pre-image $f^{-1}(K)$, or $X$ is a point. 
\end{defn}

The weight test and its stronger precursor, the coloring test, are major tools in the study of asphericity of 2-complexes. See Gersten \cite{Ger87} and Sieradski \cite{S83}. These tests can be applied to angled 2-complexes. The coloring test and one of its relative cousins is explained in detail in the next section. An important observation by Wise \cite{Wise04} states that an angled 2-complex that satisfies the coloring test has collapsing non-positive immersion. We do not know if our version of relative collapsing non-positive immersion is transitive in general. Suppose $(L,K)$ has relative collapsing non-positive immersion and $K$ has collapsing non-positive immersion. We would like to show that $L$ has collapsing non-positive immersion. So suppose $f\colon X\to L$ is a combinatorial immersion where $X$ does not have a free vertex or edge, and $X$ is not a point. Let $Y$ be the essential part of $f^{-1}(K)$. Let $Y=Y_1\cup\ldots\cup Y_n$, where the $Y_i$ are the connected components. Then
$$\chi(X)\le \chi(Y_1)+\ldots +\chi(Y_n).$$
Now $f\colon Y_i\to K$ is a combinatorial immersion, and since we assume that $K$ has collapsing non-positive immersion we have $\chi(Y_i)\le 0$, unless $Y_i$ collapses to a point. That is where the problem lies: In general we have no way to avoid collapsing connected components of $Y$. However, under certain circumstances we can control the situation. Here is one of our main results (see Theorem \ref{thm:mainapp}), vocabulary is defined in detail in the next two sections:\\

{\em Let $L$ be a standard 2-complex, $K\subseteq L$ a subcomplex all of whose 2-cells are attached along paths of exponent sum zero. Let $\bar L$ be obtained from $L$ by folding $K$ to the single edge $y$. Assume that $\bar L$ has a zero/one-angle structure that satisfies the coloring test, and $y^+$ and $y^-$ lie in different components of $lk_0(v,\bar L)$. Then $(L,K)$ has relative collapsing non-positive immersion. If in addition $K$ has collapsing non-positive immersion, then so does $L$. }\\

It was already mentioned that relative notions of combinatorial asphericity were used by the authors to show that reduced injective labeled oriented trees (LOTs) are vertex aspherical (VA) and hence aspherical. The motivation for this work is the question whether reduced injective LOTs have non-positive immersion. This is known to be true in the prime case, when the LOT is without sub-LOTs, and we answer this question affirmatively in an important but special non-prime setting. Relevant language and some history about LOTs is provided in the last section of this paper.

\section{Coloring tests}
Let $K$ be a 2-complex. If $v$ is a vertex of $K$ the {\em link at v}, $lk(v,K)$, is the boundary of a regular neighborhood of $v$. It is a graph whose edges come from the corners of 2-cells. For that reason we refer to the edges of $lk(v,K)$ as {\em corners at $v$}. If we assign numbers $\omega(c)$ to the corners $c$ of the 2-cells of $K$ we arrive at an {\em angled 2-complex}. Curvature in an angled 2-complex is defined in the following way. If $v$ is a vertex of $K$ then $\kappa(v,K)$, the curvature at $v$, is
$$\kappa(v,K)=2-\chi(lk(v,K))-\sum \omega(c_i),$$ where the sum is taken over all the corners at $v$. If $d$ is a 2-cell of $K$ then $\kappa(d,K)$, the curvature of $d$, is 
$$\kappa(d,K)=\sum \omega(c_j)-(|\partial d|-2),$$ where the sum is taken over all the corners in $d$ and $|\partial d|$ is the number of edges in the boundary of the 2-cell. The combinatorial Gauss-Bonnet Theorem states that
$$2\chi(K)=\sum_{v\in K} \kappa(v,K) + \sum_{d\in K} \kappa(d,K).$$
This was first proven by Ballmann and Buyalo \cite{BallBuy}, and later observed by McCammond and Wise \cite{McCammondWise}.
A map $X\to K$ between 2-complexes is called {\em combinatorial} if it maps open cells homeomorphically to open cells. Note that if $K$ is an angled 2-complex then the angles in the 2-cells of $K$ can be pulled back to make $X$ into an angled 2-complex. We call this angle structure on $X$ the one {\em induced} by the combinatorial map. 

An angled 2-complex where all angles are either 0 or 1 is called a {\em zero/one angled 2-complex}.
We denote by 
$lk_0(v,K)$ the subgraph of $lk(v,K)$ consisting of the vertices of $lk(v,K)$ together with the corners with angle $0$. The following coloring test is due to Sieradski \cite{S83}.

\begin{defn}(Coloring test)
Let $K$ be a zero/one-angled 2-complex. Then $K$ {\em satisfies the coloring test} if 
\begin{enumerate}
    \item the curvature of every 2-cell is $\le 0$;
    \item for every vertex $v$: If $c_1\cdots c_n$ is a reduced cycle in $lk(v,K)$, then\\ $2-\sum_{i=1}^n \omega(c_i)\le 0$. 
\end{enumerate}
\end{defn}

\begin{thm}\label{thm:npi} A zero/one-angled 2-complex $K$ satisfies the coloring test if and only if 
\begin{enumerate}
\item the curvature of every 2-cell is $\le 0$;
\item $lk_0(v,K)$ is a forest for every vertex $v$;
\item a corner with angle $1$ does not have both its vertices in a single connected component of $lk_0(v,K)$.
\end{enumerate}
\end{thm}

The proof is straightforward. 

\begin{defn}
Let $\Lambda$ be a graph. 
\begin{itemize}
    \item A cycle of edges $e_1\cdots e_n$ in $\Lambda$ is {\em reduced} if there does not exist an $e_i$ so that $e_{i+1}=\bar e_i$ ($i$ mod $n$), where $\bar e_i$ is the edge $e_i$ with reversed orientation. A cycle of edges $e_1\cdots e_n$ in $\Lambda$ is {\em homology reduced} if there does not exist a pair $e_i, e_j$ so that $e_j=\bar e_i$.
\end{itemize}    
Let $K$ be a 2-complex. 
\begin{itemize}
    \item Let $S$ be the 2-sphere. A {\em spherical diagram} is a combinatorial map $f\colon S\to K$. It is {\em reduced} if, for every vertex $v\in S$, $f$ maps $lk(v)$ to a reduced cycle. It is {\em vertex reduced} if, for every vertex $v\in S$, $f$ maps $lk(v)$ to a homology reduced cycle. 
    \item $K$ is {\em diagrammatically reducible (DR)} if there do not exist reduced spherical diagrams over $K$. $K$ is {\em vertex aspherical (VA)} if there do not exist vertex reduced spherical diagrams over $K$. 
    \end{itemize}

Let $(L,K)$ be a 2-complex pair. 
\begin{itemize}
    \item The pair is {\em relatively DR} if every reduced spherical diagram over $L$ is a diagram over $K$. It is {\em relatively VA} if every vertex reduced spherical diagram over $L$ is a diagram over $K$.
    \item A combinatorial map $f\colon X\to L$ is {\em $K$-thin} if the essential part $Y$ of $f^{-1}(K)$ has no interior vertices; that means $lk(v,Y)\ne lk(v,X)$ for every vertex $v\in X$. 
\end{itemize}
\end{defn}

\begin{thm}\label{thm:wise}
If $K$ satisfies the coloring test then it is DR and has non-positive immersion.
\end{thm}

The DR part was shown by Sieradski \cite{S83} and non-positive immersion was shown by Wise \cite{Wise04}. 

\begin{defn}\label{def:relforest}
Let $\Lambda$ be a graph and $\Lambda'$ be a subgraph, both are allowed to be disconnected. Then 
\begin{itemize}
    \item {\em $\Lambda$ is a forest relative to $\Lambda'$} if every reduced cycle in $\Lambda$ is contained in $\Lambda'$. {\em $\Lambda$ is a tree relative to $\Lambda'$} if in addition $\Lambda$ is connected.
    \item {\em $\Lambda$ is a strong forest relative to $\Lambda'$} if it is a forest relative to $\Lambda '$ and in addition, for each connected component $C$ of $\Lambda$, either $C\cap \Lambda'$ is empty or connected.
\end{itemize}
\end{defn}

\begin{figure}[htbp] 
   \centering
   \includegraphics[width=2in]{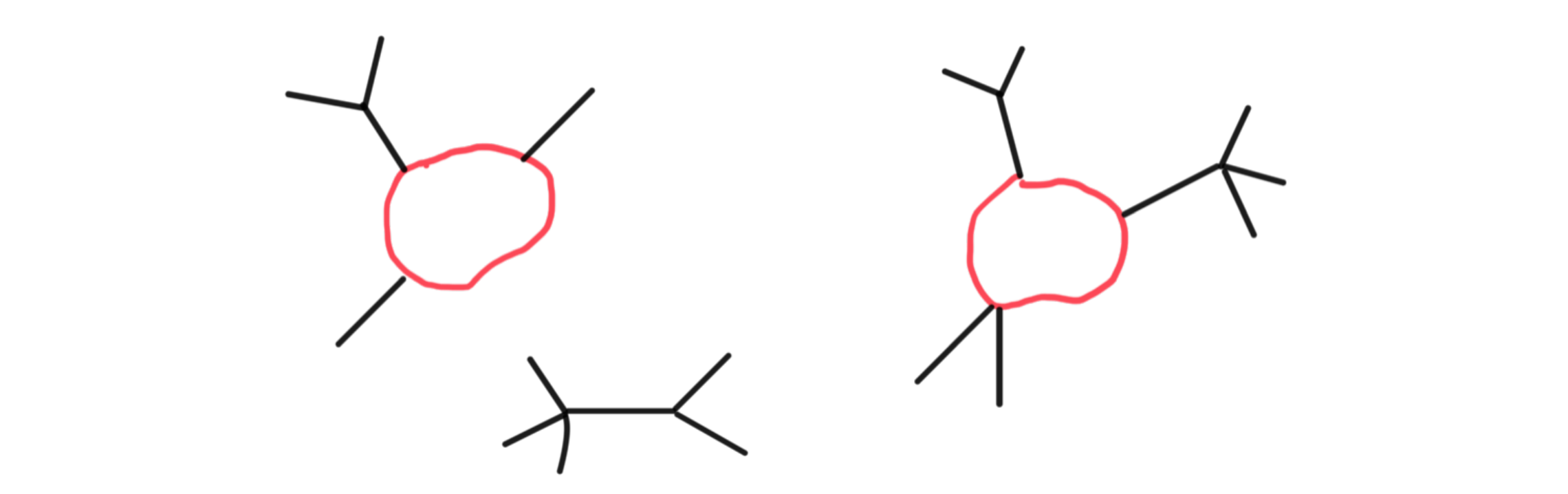} 
   \caption{A strong forest $\Lambda$ relative to $\Lambda'$ (red).}
\end{figure}

\vfill\eject

\begin{prop}\label{rem:strongrelforest} 
If $\Lambda$ is a forest relative to $\Lambda'$ and $C$ is a component of $\Lambda$ such that $C\cap \Lambda'\ne \emptyset$ then at most one component of $C\cap \Lambda'$ can not be a tree.
Let $\Lambda/\Lambda'$ be the quotient graph obtained by identifying all of $\Lambda'$ to a single vertex. Then the following conditions are equivalent:
\begin{enumerate}
    \item $\Lambda$ is a strong forest relative to $\Lambda'$; 
    \item The quotient graph $\Lambda/\Lambda'$ is a forest.
\end{enumerate}
\end{prop}

\begin{proof}
Let $C$ be a component of $\Lambda$. Suppose that $C\cap \Lambda'$ contains two components $A$ and $B$ that are not trees. Let $\alpha$ and $\beta$ be reduced closed paths in $A$ and $B$ respectively. Let $\gamma$ be a reduced path connecting $\alpha$ and $\beta$. A reduced version of $\alpha\gamma\beta\bar\gamma$ is a reduced closed path in $\Lambda$ that is not entirely contained in $\Lambda'$, contradicting the fact that $\Lambda$ is a forest relative to $\Lambda'$.

Assume that $\Lambda$ is a strong forest relative to $\Lambda'$ and let $C$ be a component of $\Lambda$. Since $C\cap \Lambda'$ is connected (or empty) and $C$ is a tree relative to $C\cap \Lambda'$, it is clear that the quotient $T_C=C/C\cap \Lambda'$ is a tree. $\Lambda/\Lambda'$ is now obtained from the forest $\bigcup T_C$ (disjoint union) by identifying single vertices from distinct $T_C$'s. This gives a forest.

Assume next that $\Lambda/\Lambda'$ is a forest. Let $C$ be a component of $\Lambda$. It is clear that $C\cap \Lambda'$ is connected or empty, because $C/C\cap \Lambda'$ is a tree. Let $q\colon C\to C/C\cap \Lambda'$ be the quotient map. Suppose that $\gamma=e_1\ldots e_n$ is a reduced closed path in $C$ not entirely contained in $C\cap\Lambda'$. Then $q(\gamma)=q(e_1)\ldots q(e_n)$ is a closed path in $C/C\cap\Lambda'$ that does contain an edge. Note that if $e_i\in \Lambda'$ then $q(e_i)$ is omitted from the sequence $q(e_1)\ldots q(e_n)$. Since $q(e_1)\ldots q(e_n)$ is a path in a tree, it contains at least two vertices of valency 1, so it contains a vertex of valency 1 that is not $C\cap \Lambda'$. 
We may assume, after cyclic reordering, that $q(e_1)$ is an edge that contains a vertex of valency 1 that is not $C\cap \Lambda'$. But then $e_1$ is an edge in $C$ that contains a vertex of valency 1, which implies that $e_1\ldots e_n$ is not reduced. In fact $e_n=\bar e_1$. A contradiction.
\end{proof}

\newpage

\begin{defn}\label{defn:relcol} (Relative coloring test)
Let $K$ be a subcomplex of a zero/one-angled 2-complex $L$. Then $(L,K)$ {\em satisfies the relative coloring test} if 
\begin{enumerate}
\item the curvature of every 2-cell $d\in L-K$ is $\le 0$;
\item for every vertex $v$: If $c_1\cdots c_n$ is a reduced cycle in $lk(v,L)$ not entirely contained in $lk(v,K)$, then $$2-\sum_{i=1}^n\omega(c_i)\le 0.$$
\end{enumerate}
\end{defn}

\begin{thm}\label{thm:relcol} A zero/one-angled pair $(L,K)$ satisfies the relative coloring test if and only if
\begin{enumerate}
\item the curvature of every 2-cell $d\in L-K$ is $\le 0$;
\item for every vertex $v$: $lk_0(v,L)$ is a forest relative to $lk_0(v,K)$;
\item if a corner $c$ with angle $1$ in a 2-cell $d\in L$ has its vertices in a single connected component of $lk_0(v,L)$, then $d\in K$ and the vertices of $c$ lie in a single component of $lk_0(v,K)$.
\end{enumerate}
\end{thm}

\begin{proof} Assume $(L,K)$ satisfies the relative coloring test. Then conditions 2 and 3 hold because otherwise one can construct a reduced cycle of corners $\lambda=c_1\ldots c_n$ not entirely contained in $lk(v,K)$ such that $\omega(\lambda)=\sum_{i=1}^n\omega(c_i)=0$ or $\omega(\lambda)=1$.

Assume that the conditions of the theorem hold.
Let $\lambda=c_1\ldots c_n$ be a reduced cycle in $lk(v,L)$ not entirely contained in $lk(v,K)$.
$\omega(\lambda)=\sum_{i=1}^n\omega(c_i)=0$ is not possible because of condition 2 in the statement of Theorem \ref{thm:relcol}. If $\omega(\lambda)=1$ then there is $1\le i\le n$ so that  $\omega(c_i)=1$ and $\omega(c_j)=0$ for $j\ne i$. Then by condition 3, $c_i$ is a corner in a 2-cell of $K$ with its vertices $p$ and $q$ in a single component of $lk_0(v,K)$. Let $\lambda'$ be a corner path in $lk_0(v,K)$ connecting $p$ to $q$. Then a reduced version of $c_i\ldots c_{i-1}\lambda'c_{i+1}\ldots c_n$ violates condition 2. Thus $\omega(\lambda)\ge 2$.
\end{proof}

The relative coloring test has implications regarding asphericity and non-positive immersion, but they are not as immediate as in the classical coloring test setting. 

\begin{thm}\label{thm:relcolasph}
Let $(L,K)$ be a 2-complex pair with a zero/one angle structure that satisfies the relative coloring test and $\kappa(d,L)\le 0$ for all 2-cells $d\in L$. If $f\colon S\to L$ is a $K$-thin spherical diagram, then it is not reduced.
\end{thm}

\begin{proof} Assume $f\colon S\to L$ is a reduced $K$-thin spherical diagram. Pull the zero/one-angle structure back to $S$. Then we have $\kappa(d,S)\le 0$ for every 2-cell $d$ of $S$. Let $v$ be a vertex in $S$. Let $c_1\cdots c_n$ be the corners that make up the link of $v$. Since $f$ is $K$-thin we have that $f(c_1)\cdots f(c_n)$ is a reduced cycle in $lk(f(v),L)$ not entirely contained in $lk(f(v), K)$. It follows that $\kappa(v,S)\le 0$. We obtain 
$$4=2\chi(S)=\sum_{v\in S}\kappa(v,S)+\sum_{d\in S}\kappa(d,S)\le 0$$
and have reached a contradiction.
\end{proof}

Under certain conditions imposed on the pair $(L,K)$ it can be shown that this result implies that $(L,K)$ is relatively DR. For more details see \cite{HarRose2020}.

\section{A strong relative coloring test}

In light of Theorem \ref{thm:relcol} we define a stronger version of the relative coloring test as follows:

\begin{defn}\label{def:strong} A zero/one-angled pair $(L,K)$ satisfies the {\em strong relative coloring test} if 
\begin{enumerate}
\item the curvature of every 2-cell $d\in L$ is $\le 0$;
\item for every vertex $v\in L$: $lk_0(v,L)$ is a strong forest relative to $lk_0(v,K)$; 
\item if a corner $c$ with angle $1$ in a 2-cell $d\in L$ has its vertices in a single connected component of $lk_0(v,L)$, then $d\in K$.
\end{enumerate}
\end{defn}

Note that if $c$ is a corner as in condition 3 with its vertices in a single component $C$ of $lk_0(v,L)$, then its vertices automatically lie in a single component of $lk_0(v,K)$ because $C\cap lk_0(v,K)$ is connected by condition 2.

Useful throughout is the following observation.

\begin{prop}\label{prop:colim}
Suppose $(L,K)$ is a zero/one-angled 2-complex pair that satisfies the (strong) relative coloring test. Let $f\colon X\to L$ be a combinatorial immersion and let $Y$ be the essential part of $f^{-1}(K)$. If we give $(X,Y)$ the angle structure induced from $(L,K)$, then $(X,Y)$ satisfies the (strong) relative coloring test.
\end{prop}

\begin{proof}
Since $f\colon X\to L$ is a combinatorial immersion it is immediate from Definition \ref{defn:relcol} that $(X,Y)$ satisfies the relative coloring test. Thus, by Theorem \ref{thm:relcol}, conditions 1 and 3 of Definition \ref{def:strong} hold for $(X,Y)$ and, for every vertex $x$, $lk_0(x,X)$ is a forest relative to $lk_0(x,Y)$.  Assume that $(L,K)$ satisfies the strong relative coloring test. In order to show that $(X,Y)$ satisfies the strong relative coloring test, the only thing left to do is show that for every vertex $x$, $lk_0(x,X)$ is not just a forest relative to $lk_0(x,Y)$, but a strong forest. Let $\Lambda=lk_0(v,L)$ and $\Lambda_K=lk_0(v,K)$. Let  $\Lambda'=lk_0(x,X)$, $\Lambda'_Y=lk_0(x,Y)$, and $\Lambda'_{f^{-1}(K)}=lk_0(x, f^{-1}(K))$, where $f(x)=v$. We can think of $\Lambda'$ is a subgraph of $\Lambda$. We have 
$$\Lambda'_Y\subseteq \Lambda'_{f^{-1}(K)}\subseteq \Lambda'\subseteq \Lambda.$$ 
Note that 
$$\Lambda'\cap \Lambda_K=\Lambda'_{f^{-1}(K)}.$$ It is easy to see that $\Lambda'$ is a strong forest relative to $\Lambda_{f^{-1}(K)}$. Indeed,
The quotient graph $\Lambda'/(\Lambda'\cap \Lambda_K)$ is a subgraph of $\Lambda/\Lambda_K$, and the latter is a forest because $\Lambda$ is a strong forest relative to $\Lambda_K$. Here we use Proposition \ref{rem:strongrelforest}. Thus $\Lambda'/(\Lambda'\cap \Lambda_K)$ is a forest and, again by Proposition \ref{rem:strongrelforest}, $\Lambda'$ is a strong forest relative to $\Lambda'_{f^{-1}(K)}$.

Since $Y$ is the essential part of $f^{-1}(K)$ we know that 
$$\Lambda'_{f^{-1}(K)}=\Lambda'_Y\ \mbox{or}\ \Lambda'_{f^{-1}(K)}=\Lambda'_Y\cup \{p_1,\ldots, p_n\},$$ 
where the $p_i$ are points not in $\Lambda_Y'$. Suppose $C$ is a component of $\Lambda'$ and $C\cap \Lambda'_Y\ne \emptyset.$ Then $C\cap \Lambda'_{f^{-1}(K)}\ne \emptyset$ and therefore $C\cap \Lambda'_{f^{-1}(K)}$ is connected because $\Lambda'$ is a strong forest relative to  $\Lambda'_{f^{-1}(K)}$. It follows that
$$C\cap \Lambda'_Y=C\cap \Lambda'_{f^{-1}(K)}$$
and so $C\cap \Lambda'_Y$ is connected.
\end{proof}

\begin{thm}\label{thm:relNPI1}
If $(L,K)$ satisfies the strong relative coloring test, then it has relative collapsing non-positive immersion.
\end{thm}

We will prove this theorem in the next section.

\bigskip Let $L$ be a 2-complex. We assume the 1-skeleton is a directed graph. If $e$ is an edge of $L$ we denote by $e^+$ a point close to the start of $e$ and by $e^-$ a point close to the end of $e$. Let $v$ be a vertex of $L$ and let $E(v)$ be the set of edges starting or ending at $v$. We denote by $lk^+(v,L)$ and $lk^-(v,L)$ the subgraphs of $lk(v,L)$ spanned by the $e^+$, and the $e^-$, $e\in E(v)$, respectively. Corners that are neither positive nor negative (i.e. not in $lk^+(v,L)\cup lk^-(v,L)$) are called mixed corners. 

There is a {\em standard zero/one-angle assignment} for $L$: Give positive and negative corners angle $0$, and mixed corners angle $1$. Note that in this case we have
$$lk_0(v,L)=lk^+(v,L)\cup lk^-(v,L),\ lk^+(v,L)\cap lk^-(v,L)=\emptyset.$$ The second statement always holds and is independent of angle assignments.

Assume $L$ is a 2-complex and $K$ is a subcomplex all of whose 2-cells are attached along closed paths of exponent sum zero. We attach 2-cells to every closed path of exponent sum zero in $K$, reduced or not,  to obtain a larger 2-complex $\hat K$. Let $\hat L$ be $L$ together with these added 2-cells. We call $(\hat L,\hat K)$ the {\em maximal expansion of $(L,K)$}.  Note that $lk^+(v, \hat K)$ contains a full subgraph on its vertices, and so does $lk^-(v,\hat K)$. In particular both $lk^+(v, \hat K)$ and $lk^-(v, \hat K)$ are connected.

\begin{thm}\label{thm:relNPI}
Let $L$ is a 2-complex and $K$ is a subcomplex all of whose 2-cells are attached along loops of exponent sum zero. Assume 
\begin{enumerate}
    \item $(\hat L,\hat K)$ carries a zero/one-angle structure which satisfies the strong relative coloring test;
    \item the angle structure on $\hat K$ is standard;
    \item $K$ has collapsing non-positive immersion.
\end{enumerate}
Then $L$ has collapsing non-positive immersion.
\end{thm}

Since the maximal expansion $(\hat L,\hat K)$ satisfies the strong relative coloring test it has relative collapsing non-positive immersion by Theorem \ref{thm:relNPI1}. Since this property is hereditary by Proposition \ref{prop:colim}, it follows that $(L,K)$ has collapsing non-positive immersion. 
Theorem \ref{thm:relNPI} established a transitivity property for collapsing non-positive immersion in a special setting.
We will prove this theorem in a later section. It will involve a process of \lq\lq thinning\rq\rq\  a given combinatorial immersion.\\

\begin{figure}[ht]\centering
\begin{tikzpicture}[scale=1.3]

\begin{scope}[decoration={markings, mark=at position 0.5 with {\arrow{>}}}]
\draw [postaction={decorate}] (0,0) -- (1,0) node[midway, below]{$y_2$};
\draw [postaction={decorate}] (1,0) -- (1,1) node[midway, right]{$y_1$};
\draw [postaction={decorate}] (0,0) -- (0,1) node[midway, left]{$y_3$};
\draw [postaction={decorate}] (0,1) -- (1,1) node[midway, above]{$y_2$};
\draw [postaction={decorate}] (2,0) -- (2.8,0.5) node[midway, below]{$y$};
\draw [postaction={decorate}] (2.8,0.5) -- (3.1,1.5) node[midway, right]{$y$};
\draw [postaction={decorate}] (2,0) -- (2.3,1) node[midway, left]{$y$};
\draw [postaction={decorate}] (2.3,1) -- (3.1,1.5) node[midway, above]{$y$};
\draw [postaction={decorate}] (4,0) -- (4.4,1) node[midway, right]{$y$};
\end{scope}
\draw[blue] (2.3,0.2) arc (30:45:1cm);
\draw[blue] (2.8,1.3) arc (210:225:1cm);
\fill (4,0) circle (1.6pt);
\fill (4.4,1) circle (1.6pt);
\fill[blue] (4.08,0.2) circle (1.3pt); \node[right,blue] at (4.08,0.2) {$y^+$};
\fill[blue] (4.3,0.8) circle (1.3pt); \node[left,blue] at (4.3,0.8) {$y^-$};
\end{tikzpicture}
\caption{\label{afoldy} A 2-cell in $K$ is folded to the edge $y$.}
\end{figure}
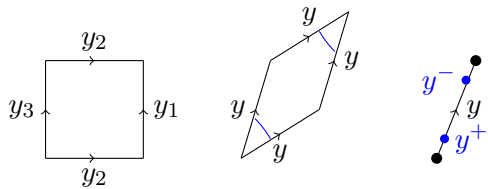

Assume $L$ is a standard 2-complex with single vertex $v$, and $K$ is a subcomplex. Assume that all 2-cells of $K$ are attached along paths of exponent sum zero. Let $y$ be an edge of $K$. Notice that we can ``fold" $K$ onto $y$ by mapping every edge of $K$ to $y$ and extending this map to the 2-cells (see Figure \ref{afoldy}). Let $\bar L$ be the quotient of $L$ obtained by folding $K\subseteq L$ to $y$. 

\begin{example} Let 
$$P_0=\langle y_1,y_2\ |\  y_1y_2y_1^{-1}y_2^{-1}\rangle \subseteq P=\langle x, y_1,y_2\ |\  y_1y_2y_1^{-1}y_2^{-1}, xy_1xy_2 \rangle.$$ Let $K\subseteq L$ be the presentation complexes of $P_0$ and $P$, respectively. We can fold $K$ to $y$ and obtain $\bar L$, the presentation complex for $\langle x, y\ | xyxy \rangle$.
\end{example}

Here is a key observation. $lk(v,\bar L)$ is obtained from $lk(v,L)$ by the following process: 
\begin{enumerate}
    \item Remove all mixed corners in $lk(v,K)$ from $lk(v,L)$.
    \item Identify all of $lk^+(v,K)$ with $y^+$, and all of $lk^-(v,K)$ with $y^-$.
\end{enumerate}

\noindent Here is our main application of Theorem \ref{thm:relNPI}. It was already mentioned in the introduction.

\begin{thm}\label{thm:mainapp}
Let $L$ be a standard 2-complex, $K\subseteq L$ a subcomplex all of whose 2-cells are attached along paths of exponent sum zero. Let $\bar L$ be obtained from $L$ by folding $K$ to the single edge $y$. Assume that $\bar L$ has a zero/one-angle structure that satisfies the coloring test, and $y^+$ and $y^-$ lie in different components of $lk_0(v,\bar L)$. Then $(L,K)$ has relative collapsing non-positive immersion. If in addition $K$ has collapsing non-positive immersion, then so does $L$. 
\end{thm}

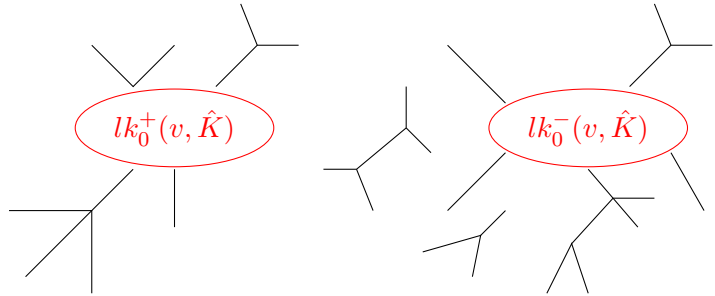
\begin{figure}[ht]\centering
\begin{tikzpicture}[scale=1.1]
\node[ellipse, draw, red] (e) at (0,0) {$lk_0^+(v,\hat K)$};
\node[ellipse, draw, red] (e) at (5,0) {$lk_0^-(v,\hat K)$};
\draw (-0.5,-0.5) -- (-1.8,-1.8);
\draw (-1,-2) -- (-1,-1) -- (-2,-1);
\draw (0,-0.5) -- (0,-1.2);
\draw (0,1) -- (-0.5, 0.5) -- (-1, 1);
\draw (0.5, 0.5) -- (1,1) -- (1.5,1);
\draw (1,1) -- (0.8,1.5);
\draw (1.8,-0.5) -- (2.2,-0.5) -- (2.8,0) -- (2.8,0.5);
\draw (2.2,-0.5) -- (2.4, -1);
\draw (2.8,0) -- (3.1,-0.3);
\draw (3,-1.5) -- (3.7,-1.3) -- (3.6, -1.8);
\draw (3.7,-1.3) -- (4,-1);
\draw (5,-0.5) -- (5.6,-1.2);
\draw (4.8,-1.4) -- (5.3,-0.85) -- (5.8,-0.85);
\draw (4.5,-2) -- (4.8,-1.4) -- (5,-2);
\draw (6, -0.3) -- (6.4,-1);
\draw (4, -0.3) -- (3.3, -1);
\draw (5.5, 0.5) -- (6,1) -- (6.5,1);
\draw (6,1) -- (5.8,1.5);
\draw (4, 0.3) -- (3.3, 1);
\end{tikzpicture}
\caption{\label{fig:linkstrong}$lk_0(v,\hat L)$}
\end{figure}

\begin{proof}
Give $\bar L$ a zero/one-angle structure that satisfies the coloring test. Assign to the maximal expansion $\hat K$ the standard zero/one-angle structure, and for a 2-cell in 
$\hat L$ not in $\hat K$ the angle structure pulled back from $\bar L$.
We will show that $(\hat L,\hat K)$ satisfies the strong relative coloring test and therefore the conditions of Theorem \ref{thm:relNPI} hold. Note first that we also have a map $\hat L\to \bar L$, folding all of $\hat K$ to the single edge $y$.\\

\noindent 1. If $d$ is a 2-cell in $\hat K$ then $\kappa(d,\hat K)\le 0$. This is because 2-cells in $\hat K$ are attached along paths of exponent sum zero and $\hat K$ carries the standard zero/one-angle structure. If $d$ is a 2-cell in $\hat L$ not in $\hat K$, then $\kappa(d,\hat L)\le 0$ because $\kappa(\bar d,\bar L)\le 0$, where $\bar d$ is the image of $d$ under the folding map. \\

\noindent 2. We have a quotient map $$q\colon lk_0(v,\hat L) \to lk_0(v,\bar L),$$
where all of $lk^+(v,\hat K)$ gets mapped to $y^+$, and all of $lk^-(v,\hat K)$ gets mapped to $y^-$. Note that $lk_0(v,\bar L)$ is a forest because $\bar L$ satisfies the coloring test.
Since $y^+$ and $y^-$ lie in different components of $lk_0(v, \bar L)$, $lk^+(v,\hat K)$ and $lk^-(v,\hat K)$ lie in different components of $lk_0(v,\hat L)$. The situation is shown in Figure \ref{fig:linkstrong}. Let $C^+$ be the component of $lk_0(v,\hat L)$ that contains $lk^+(v,\hat K)$. Since $q(C^+)$ is a tree, it follows that $C^+$ is a tree relative to $lk^+(v,\hat K)$. This follows from Proposition \ref{rem:strongrelforest}. Similarly, $C^-$ is a tree relative to $lk^-(v,\hat K)$. It follows that $lk_0(v,\hat L)$ is a forest relative to $lk_0(v,\hat K)=lk^+(v,\hat K)\cup lk^-(v,\hat K)$. Furthermore, if $C$ is any component of $lk_0(v,\hat L)$, then either $C\cap lk_0(v,\hat K)=\emptyset$, or $C\cap lk_0(v,\hat K)=lk_0^+(v,\hat K)$, or $C\cap lk_0(v,\hat K)=lk_0^-(v,\hat K)$. So $C\cap lk_0(v,\hat K)$ is either empty or connected.\\ 

\noindent 3. Let $c$ be a corner of angle 1. If $c$ comes from a 2-cell in $\hat K$, then $c$ connects $lk^+(v, \hat K)$ to $lk^-(v, \hat K)$, which lie in different components of $lk_0(v,\hat L)$. Assume $c$ comes from a 2-cell $d$ in $\hat L$ not in $\hat K$. The folding map $\hat L\to \bar L$ sends $c$ to $\bar c$ in $\bar d$. $\bar c$ is a corner of angle 1 and hence connects two distinct components of $lk_0(v,\bar L)$, because $\bar L$ satisfies the coloring test. Therefore $c$ connects distinct components of $lk_0(v,\hat L)$.
\end{proof}

\pagebreak

\begin{cor}\label{cstnpi}
 Let $L$ be a standard 2-complex, $K\subseteq L$ a subcomplex all of whose 2-cells are attached along paths of exponent sum zero. Let $\bar L$ be obtained from $L$ by folding $K$ to the single edge $y$. Assume that $\bar L$ satisfies the coloring test with the standard zero/one-angle structure. Then $(L,K)$ has relative collapsing non-positive immersion. If in addition $K$ has collapsing non-positive immersion, then so does $L$. 
\end{cor}

\begin{proof}
Note that in this setting $lk_0(v,\bar L)=lk^+(v,\bar L)\cup lk^-(v,\bar L)$. We have $y^+\in lk^+(v,\bar L)$, $y^-\in lk^-(v,\bar L)$, and $lk^+(v,\bar L)\cap lk_0^+(v,\bar L)=\emptyset$. Thus $y^+$ and $y^-$ lie in different components of $lk_0(v, \bar L)$.
\end{proof}

\begin{example}
Consider a presentation $P_0=\langle {\bf y}\ |\ {\bf r}\rangle$, where every $r\in {\bf r}$ has exponent sum zero. Let $P=\langle x,  {\bf y}\ |\ {\bf r}, xux^{-1}v^{-1} \rangle$, where $u$ and $v$ are words in ${\bf y}^{\pm 1}$ of equal exponent sum $k>0$. The presentation $P$ typically arises in the context of HNN-extensions of groups, where the associated subgroup is $\mathbb Z$. Let $L$ be the presentation complex of $P$, and let $K$ be the subcomplex coming from $P_0$. Choose $y\in {\bf y}$. We fold $K$ onto $y$ and obtain a quotient $L\to \bar L$. Note that $\bar L$ is the standard 2-complex constructed from $\langle x,y\ |\ xy^kx^{-1}y^{-k}\rangle$. It is easy to check that $\bar L$ satisfies the coloring test with the standard zero/one-angle structure. It follows by Corollary \ref{cstnpi} that $(L,K)$ has relative non-positive immersion.
\end{example}

\section{Proof of Theorem \ref{thm:relNPI1}}
\begin{prop}\label{prop:heredi}
Assume that $(L,K)$ satisfies the strong relative coloring test and that $L$ has no free vertex or edge.  Then $\chi(L)\le \chi(K)$. Furthermore, if there exist vertices $v_1,\ldots, v_n\in K$ such that $\kappa(v_i,L)<\kappa(v_i,K)$, then $\chi(L)\le \chi(K)-n$.
\end{prop}

\begin{proof}
Let $v$ be a vertex of $L$. We have 
$$\kappa(v,L)=2-\chi(lk(v,L))-\sum \omega(c_i)=2-\chi(lk_0(v,L)),$$
where the sum is taken over all corners in the link of $v\in L$.
Assume first that $v$ is a vertex of $L$ that is not in $K$. Then every component of $lk_0(v,L)$ is a tree. We claim that there has to be more than one component. $lk(v,L)$ can not have vertices of valency 0 or 1 because $L$ does not contain a free vertex or edge. Thus a vertex of valency 0 or 1 in $lk_0(v,L)$ is the vertex of a corner with angle 1. That corner connects distinct components. Therefore $lk_0(v,L)$ has at least two components. It follows that $\chi(lk_0(v,L))\ge 2$ and therefore
$$\kappa(v,L)=2-\chi(lk_0(v,L))\le 0.$$

Next assume that $v$ is a vertex of $K$. Assume $lk_0(v,L)=C_1\cup\dots\cup C_n$, where the $C_i$ are the connected components. Assume first that $C_i\cap lk_0(v,K)=\emptyset$ for $1\le i\le k$ and $C_i\cap lk_0(v,K)\ne \emptyset$ for $k+1\le i \le n$ ($k=0$ and $k=n$ are both possible). Then $C_i$ is a tree for $1\le i\le k$, and $C_i$ is a tree relative to the connected subgraph $C_i\cap lk_0(v,K)$ for $k+1\le i \le n$. Thus $C_i$ is homotopically equivalent to $C_i\cap lk_0(v,K)$ for $k+1\le i \le n$. We have 
\begin{equation}\label{eq:first}
    \chi(lk_0(v,L))=k+\sum_{i=k+1}^n \chi(C_i\cap lk_0(v,K))\ge k+\chi(lk_0(v,K)).
\end{equation}
Thus
\begin{equation}\label{eq:second}
\kappa(v,L)=2-\chi(lk_0(v,L))\le 2-(\chi(lk_0(v,K))+k)=\kappa(v,K)-k.
\end{equation}
Using $\kappa(v,L)\le 0$, $v\not\in K$, and $\kappa(d,L)\le 0$, $d\not\in K$, we obtain 
\begin{equation}\label{eq:third}
2\chi(L)=\sum_{v\in L} \kappa(v, L)+\sum_{d\in L} \kappa(d,L)\le \sum_{v\in K} \kappa(v,K)+\sum_{d\in K} \kappa(d,K)=2\chi(K).
\end{equation}
Note that if there exist vertices  $v_1,\ldots v_n\in K$ such that $\kappa(v_i,L)<\kappa(v_i,K)$, then $\chi(L) \le \chi(K)-n$ follows immediately from the last Equation \ref{eq:third}. 
\end{proof}

{\bf Proof of Theorem \ref{thm:relNPI1}}: 
Let $f\colon X\to L$ be a combinatorial immersion. Assume that $X$ is finite, connected, has no free vertex or edge, and that $X$ is not a single point. Let $Y$ be the essential part of $f^{-1}(K)$. Note that $(X,Y)$ satisfies the strong relative coloring test by Proposition \ref{prop:colim}. Therefore $\chi(X)\le \chi(Y)$ by Proposition \ref{prop:heredi}. \qed

\section{A collection of lemmas needed for Theorem \ref{thm:relNPI}}

Let $(L,K)$ be a 2-complex pair. An {\em interior point $p$ of $K$} is a point so that $lk(p,L)=lk(p,K)$. Define $K^{\circ}$ to be the set of interior points and $\partial K = K-K^{\circ}$. Note that $\partial K=(L-K^{\circ})\cap K$ is a subgraph  of the 1-skeleton of $K$. 
If all of $K$ is essential, that is every edge is part of the attaching path of some 2-cell, then every edge in $\partial K$ is in the boundary of a 2-cell of $K$ and a 2-cell of $L-K^{\circ}$. 

\pagebreak

\begin{lemma}\label{lem:val1}
Suppose $(L,K)$ satisfies the strong coloring test, $L$ has no free vertex or edge. We assume that all of $K$ is essential. Let $v$ be a vertex of $\partial K$. Assume
\begin{enumerate}
    \item $v$ is an isolated vertex (valency 0) in $\partial K$, or
    \item $v$ has valency 1 in $\partial K$, or  
    \item $lk_0(v,K)$ is connected.
\end{enumerate}
Then $lk_0(v,L)$ contains a connected component $C$ so that $C\cap lk_0(v,K)=\emptyset$. Thus, in all three cases $\kappa(v,L)<\kappa(v,K)$, and therefore we have $\chi(L)<\chi(K)$.

\end{lemma}

\begin{proof}
Let $C$ be a component of $lk_0(v,L)$ such that $C\cap lk_0(v,K)=\emptyset$. Then $\kappa(v,L)<\kappa(v,K)$. This was shown in the second paragraph of the proof of Proposition \ref{prop:colim}, see Equation \ref{eq:second}.
That this implies that $\chi(L)<\chi(K)$ is part of the same Proposition \ref{prop:colim}.\\

\noindent 1. Here we have $lk(v,L)=lk(v,K)\cup lk(v, L-K^\circ)$ and  $lk(v,K)\cap lk(v, L-K^\circ)=\emptyset$. Since $v\in \partial K$ the link  $lk(v, L-K^\circ)\ne \emptyset$. Let $p$ be a vertex in $lk(v,L-K^{\circ})$ and let $C$ be the component of $lk_0(v,L)$ that contains $p$. Then $C\subseteq lk(v, L-K^\circ)$ and therefore $C\cap lk_0(v,K)=\emptyset$.\\

\noindent 2. Let $e$ be the unique edge in $\partial K$ which contains $v$. Let $\{ p\} =lk(v,L)\cap e$. We have $$lk(v,L)=lk(v,K)\cup lk(v, L-K^\circ) \ \mbox{and} \ lk(v,K)\cap lk(v, L-K^\circ)=\{ p \}$$ and this implies  
$$lk_0(v,L)=lk_0(v,K)\cup lk_0(v, L-K^\circ) \ \mbox{and} \ lk_0(v,K)\cap lk_0(v, L-K^\circ)=\{ p \}.$$ 
Thus, if $C$ is a component of $lk_0(v,L)$ so that $p\ne C$ and $C\cap lk(v, L-K^\circ)\ne \emptyset$, then $C\cap lk(v,K)=\emptyset$.

Since $e\in \partial K$, there is a corner in $lk(v,L-K^{\circ})$ with endpoints $p$ and $q$. Note that $p\ne q$ because $(L,K)$ satisfies the strong coloring test. Because $v$ has valency 1 in $\partial K$ we have $q\not\in lk(v,K)$. Let $C_p$ and $C_q$ be the components of $lk_0(v,L)$ that contain $p$ and $q$, respectively. See Figure \ref{fig:lemmalink}.

Assume first that that $C_p\ne C_q$.
We have $p\not\in C_q$ and 
$C_q\cap lk(v,L-K^{\circ})\ne \emptyset$ because the intersection contains $q$. Thus $C_q\cap lk_0(v,K)=\emptyset$.

Assume next that $C_p=C_q$.
Since $C_p$ contains $q\not\in lk(v,K)$, $C_p$ is not contained in $lk_0(v,K)$. Since $C_p$ is a tree relative to $C_p\cap lk_0(v,K)$, it contains a vertex $r\not\in lk_0(v,K)$ of valency 1. Since $L$ does not have a free edge there is a corner $c$ of angle 1 with vertex $r$. Since $c\not\in lk(v,K)$ (because $r\not\in lk(v,K)$), $c$ connects $C_p$ to a different component $C$ of $lk_0(v,L)$: $C_p\cap C=\emptyset$. Here we used condition 3 in of Definition \ref{def:strong}. In particular $p\not\in C$. 
Let $r$ and $s$ be the vertices of the corner $c$. Since $c\in lk(v,L-K^{\circ})$, $s\in lk(v,L-K^{\circ})$. We have $p\not\in C$ and $C\cap lk(v,L-K^{\circ})\ne \emptyset$ because the intersection contains $s$. Thus $C\cap lk_0(v,K)=\emptyset$.\\

\begin{figure}
\centering
\begin{subfigure}{.4\textwidth}
  \centering
  \includegraphics[width=1.1\linewidth]{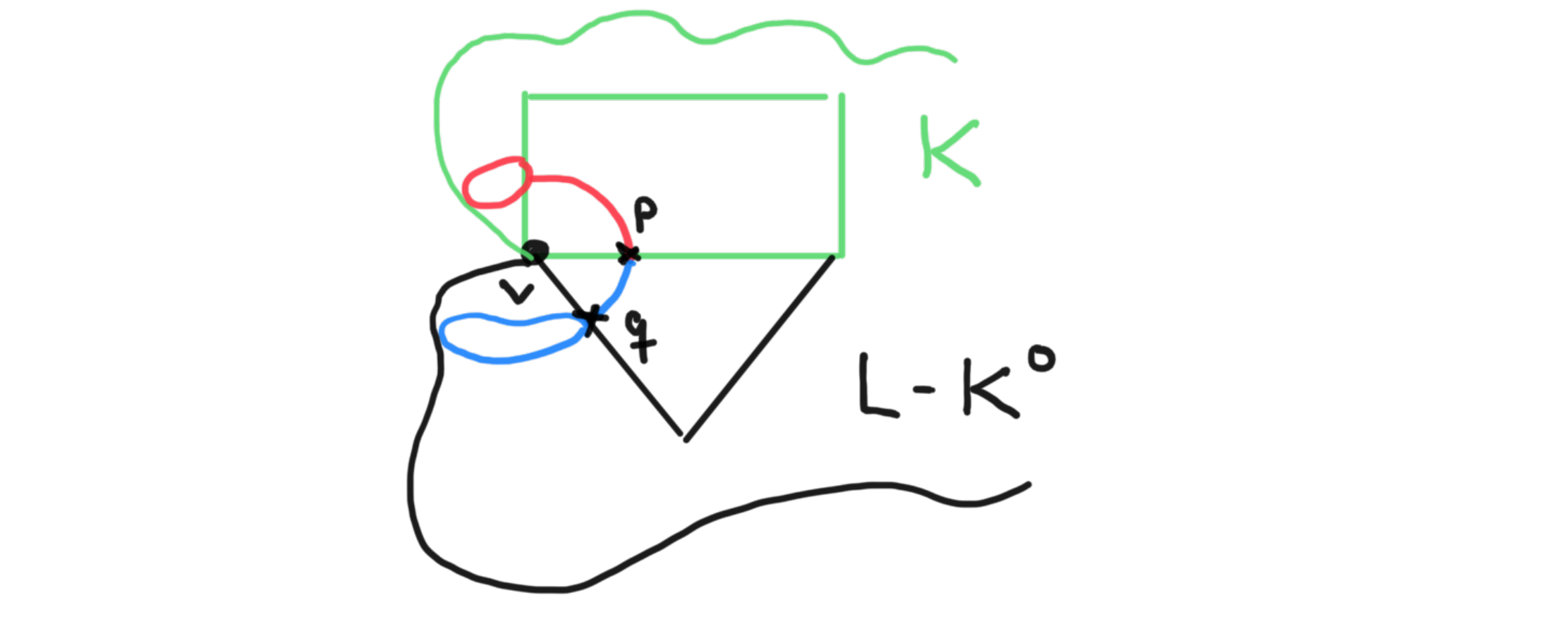}
  \label{fig:sub1}
\end{subfigure}%
\begin{subfigure}{.4\textwidth}
  \centering
  \includegraphics[width=1.1\linewidth]{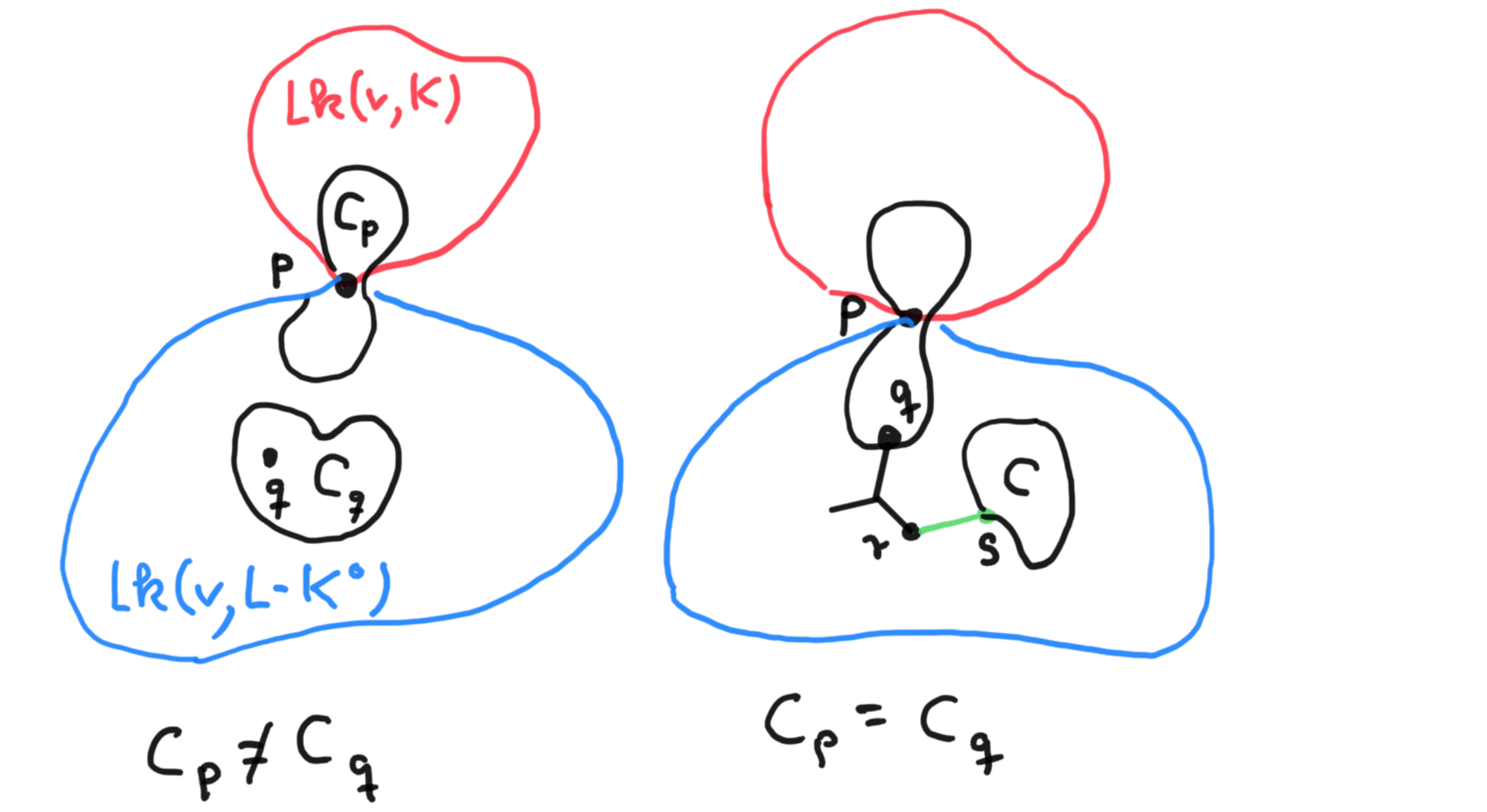}
  \label{fig:sub2}
\end{subfigure}
\caption{On the left we see a vertex $v$ of valency 1 in $\partial K$. In this case $lk(v,L-K^{\circ}\cap lk(v,K)=\{ p\}$. On the right there is the construction of a component $C$ of $lk_0(v,L)$ such that $C\cap lk_(v,K)=\emptyset$. This construction distinguishes two cases: $C_p\ne C_q$ and $C_p=C_q$.}
\label{fig:lemmalink}
\end{figure}

\noindent 3. Let $C$ be the component of $lk_0(v,L)$ that contains $lk_0(v,K)$. If there is another component $C'$ then $C'\cap lk_0(v,K)=\emptyset$ and we are done. So suppose $C$ is the only component $C=lk_0(v,L)$. 
Then $lk(v,L-K^{\circ})$ can not contain corners of angle 1, because such connect distinct components. Thus $lk(v,L)=C\cup lk(v,K)$. Since $C$ is a tree relative to $lk_0(v,K)$, it follows that $lk(v,L)$ is a tree relative to $lk(v,K)$. Since $v\in \partial K$, $lk(v,K)\ne lk(v,L)$. Therefore
$lk(v,L)$ contains a vertex of valency 1, contradicting the fact that $L$ does not have a free edge.
\end{proof}

\begin{lemma}\label{lem:attach}
Suppose $\Lambda$ is a finite connected graph that is not a tree. Then 2-cells can be attached to obtain a 2-complex $K$ so that:
\begin{enumerate}
    \item $K$ is collapsible;
    \item every edge is in the boundary of a 2-cell;
    \item $lk(v,K)$ is connected for every vertex $v$.
\end{enumerate}
\end{lemma}

\begin{proof}
We will do induction on the number of edges in $\Lambda$. If there is only one edge then $\Lambda$ is a circle and the statement is true. 

Choose a maximal tree $T$ in $\Lambda$. Remove the interior of an edge $e\not\in T$ from $\Lambda$ to obtain $\Lambda'$. If $\Lambda'$ is a tree then attach a 2-cell along the path $e\gamma$ to $\Lambda $, where $\gamma$ is the path that runs around the tree, from the endpoint of $e$ to its starting point, visiting all vertices. See Figure \ref{fig:2cell}. It is easy to check that this produces a 2-complex with the desired properties.

\begin{figure}[htbp] 
   \centering
   \includegraphics[width=2in]{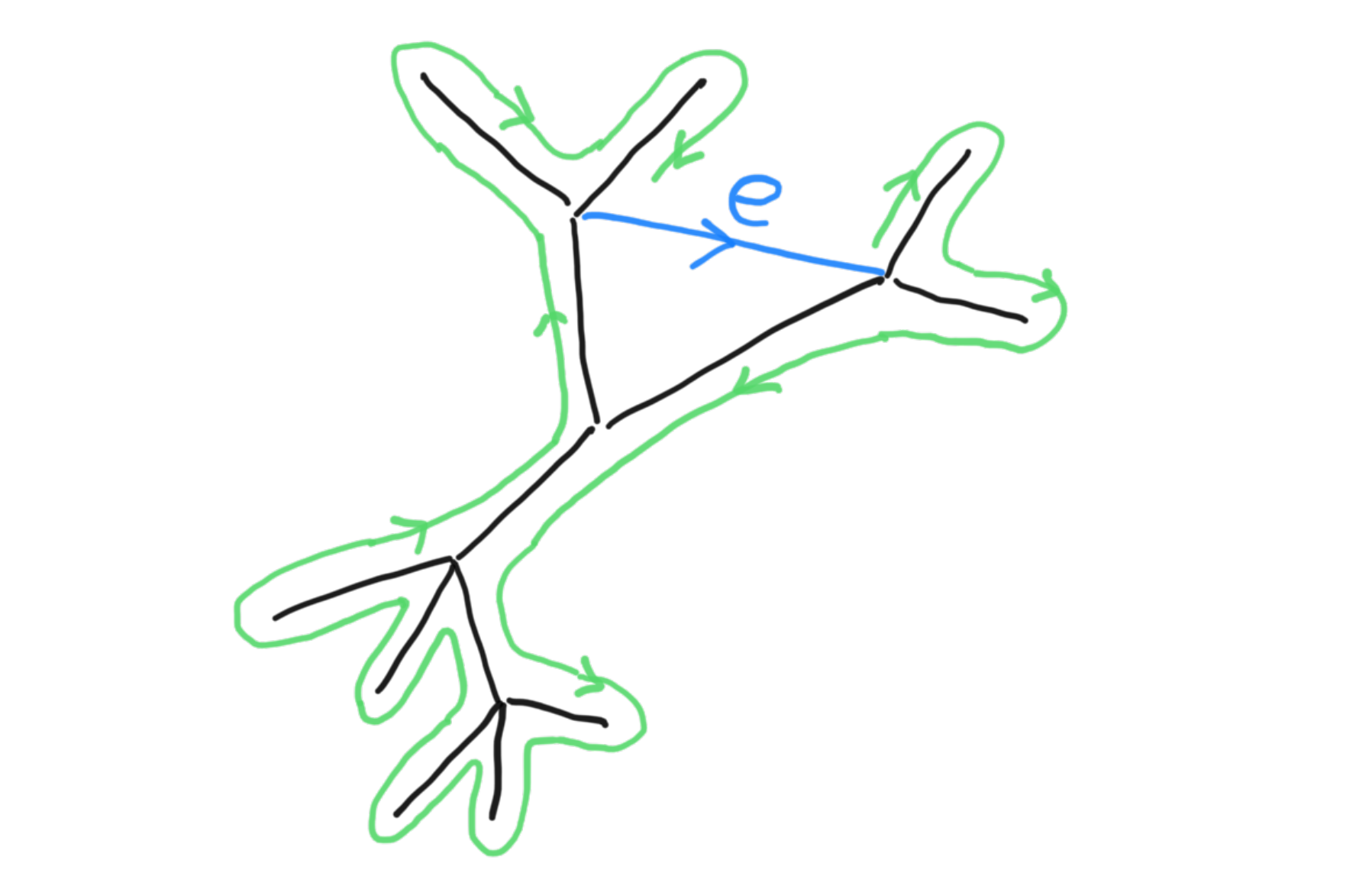} 
   \caption{\label{fig:2cell}The 2-cell is attached along $e\gamma$, where $\gamma$ goes around a maximal tree covering every edge. Note that $\gamma$ need not be reduced.}
\end{figure}

Now assume $\Lambda'$ is not a tree. Then by induction we can attach 2-cells to $\Lambda'$ and produce a 2-complex $K'$ with the desired properties. Let $p$ and $q$ be the starting and the ending vertex of $e$ ($p=q$ is possible) respectively. We attach a 2-cell $d$ to $K'$ along the path $e\gamma$, where $\gamma$ is a path in $T$ from $q$ to $p$ that covers all of $T$. See Figure \ref{fig:2cell}. This gives a 2-complex $K$. Condition 1 is clearly satisfied. $K$ is collapsible because the last 2-cell we attached has $e$ as a free edge and $K'$ is collapsible. Lets check the links. If $v$ is a vertex different from both $p$ and $q$, then adding the 2-cell $d$ adds edges to $lk(v,K')$, but not any new vertices. Since $lk(v,K')$ is connected, so is $lk(v,K)$. Adding $d$ introduces a new vertex $e^+$ to the link at $p$, but also an edge connecting $e^+$ to $lk(p,K')$. Again, since $lk(p,K')$ is connected, so is $lk(p,K)$. The same argument also shows that $lk(q,K)$ is connected.
\end{proof}

A vertex $v$ in a 2-complex $K$ is called a {\em sink (source)} if all edges in $K$ that contain $v$ point towards (away from) $v$.

\begin{lemma}\label{lem:sink}
Let $K$ be a connected finite 2-complex all of whose 2-cells are attached along paths of exponent sum zero. Assume $\pi_1(K)=1$. Then every closed path in $K$ has exponent sum zero. Furthermore, $K$ has sink and source vertices.
\end{lemma}

\begin{proof}
Because $\pi_1(K)=1$ every closed path is, up to free reductions and expansions, of the form 
$$\prod \gamma_i(\partial d_i)^{\pm 1}\gamma_i^{-1}$$
where the $d_i$ are not necessarily distinct 2-cells. Since $\partial d_i$ has exponent sum zero, the first statement follows.

Fix a vertex $v_0$. We define a map $h\colon V(K)\to \mathbb Z$, where $V(K)$ is the set of vertices of $K$, in the following way: If $v$ is a vertex of $K$ and $\gamma$ is a path from $v_0$ to $v$, let $h(v)$ be the exponent sum of $\gamma$. This is well defined by the first statement just shown. Let $v_m$ and $v_M$ be vertices where $h$ is minimal and maximal, respectively. Then $v_m$ is a source and $v_M$ is a sink.
\end{proof}

\section{Proof of Theorem \ref{thm:relNPI}}

Let $(L,K)$ be a zero/one-angled 2-complex pair so that the conditions of Theorem \ref{thm:relNPI} hold: 
\begin{enumerate}
    \item $L$ is a 2-complex and $K$ is a subcomplex all of whose 2-cells are attached along loops of exponent sum zero.
    \item The expansion $(\hat L,\hat K)$ carries a zero/one-angle structure which is standard on $\hat K$ and which satisfies the strong relative coloring test.
    \item $K$ has collapsing non-positive immersion.
\end{enumerate}
We will show that $L$ has collapsing non-positive immersion.\\

Let $f\colon X\to L$ be a combinatorial immersion. We assume $X$ does not contain free edges or vertices and is not a point. We will show that $\chi(X)\le 0$. Let $Y$ be the essential part of $f^{-1}(K)$. Give $(X,Y)$ the zero/one-angle structure induced from $(L,K)\subseteq (\hat L,\hat K)$. Since $(\hat L,\hat K)$ satisfies the strong relative coloring test, Proposition \ref{prop:colim} implies that both $(L,K)$ and $(X,Y)$ do as well.
Also note that because the immersion $f$ is combinatorial and cells in $K$ are attached along paths of exponent sum zero, cells of $Y$ are attached along paths of exponent sum zero. And since the angle structure on $K$ is standard, the angle structure on $Y$ is standard as well.  Let $Y=Y_1\cup\dots \cup Y_n$ be the disjoint union of connected components. Each $f\colon Y_i\to K$ is an immersion. Since we assume that $K$ has collapsing non-positive immersion it follows that $\chi(Y_i)\le 0$ or $Y_i$ is collapsible. Thus, if no $Y_i$ is collapsible, then, by Proposition \ref{prop:heredi}, 
$$\chi(X)\le \chi(Y)=\sum_{i=1}^n\chi(Y_i)\le 0$$ and we are done. Assume
$$Y=Y_1\cup\ldots\cup Y_m\cup\ldots\cup Y_n \ \mbox{(connected components)}$$
where the $Y_i$, $1\le i\le m$, are collapsible, and the $Y_i$, $m+1\le i\le n$, are not. The idea now is to find vertices $v_i\in Y_i$, $1\le i\le m$ so that $\kappa(v_i,X)<\kappa(v_i,Y)$. Then, again by Proposition \ref{prop:heredi}, we have
$$\chi(X)\le \chi(Y)-m=\sum_{i=1}^m\chi(Y_i)-m+\sum_{i=m+1}^n\chi(Y_i).$$
Since $\chi(Y_i)=1$, $1\le i\le m$, and $\chi(Y_i)\le 0$, $m+1\le i\le n$, we obtain the desired result $\chi(X)\le 0$. Lemma \ref{lem:val1} will be crucial in locating the vertices $v_i$. There is one case however where this lemma is of little help: If for some $1\le i\le m$ the graph $\partial Y_i$ does not contain a vertex $v$ of valency 0 or 1, or where $lk_0(v,Y)$ is not connected. In order to avoid this case we will carefully replace the combinatorial immersion $f\colon X\to L$ by a more convenient one $f'\colon X'\to \hat L$ such that $\chi(X)\le \chi(X')$. 
Assume 
$$Y=Y_1\cup\ldots\cup Y_k\cup \ldots\cup Y_m\cup\ldots\cup Y_n \ \mbox{(connected components)},$$
where $Y_i$ is collapsible for $1\le i\le m$, $\partial Y_i$, $1\le i\le k$, does not contain a vertex $v$ of valency 0 or 1, but $\partial Y_j$, $k+1\le j\le m$, does.
Let $$Y_{col}=Y_1\cup\ldots\cup Y_k$$ and 
$$X_{rest}=X-Y_{col}^{\circ}.$$ 
Assume $$\partial Y_{col}=\Delta_1\cup\ldots\cup\Delta_p \ \mbox{(connected components)}.$$
Note that none of the $\Delta_i$, $1\le i\le p$, is a tree, because we assumed that $\partial Y_i$, $1\le i\le k$, does not contain a vertex or valency 0 or 1.  $X_{rest}$ might not be connected and the $\Delta_i$ are distributed over the various connected components of $X_{rest}$. 

We attach 2-cells to each $\Delta_i$, $1\le i\le p$, to obtain 2-complexes $Z_i$ that satisfy the conditions of Lemma \ref{lem:attach}. Let 
$$Z=Z_1\cup\ldots\cup Z_p$$ and 
$X'=X_{rest}\cup Z$. Note that $Z$ is essential (meaning every edge is part of the attaching path of a 2-cell in $Z_i$). Thus every edge in $Z$ belongs to a 2-cell of $Z$ and a 2-cell of $X_{rest}$. This implies that $Z^{\circ}=\mbox{\{open 2-cells in $Z$\}}$, that the 1-skeleton of $Z$ is $\partial Z$, and that $X'$ does not have a free vertex or a free edge. 

We will next extend the immersion $f|_{X_{rest}}\colon X_{rest}\to L$ to an immersion $f'\colon X'\to \hat L$. We have already mentioned that cells in $Y$ are attached along paths of exponent sum zero. Since the components of $Y_{col}$ are collapsible, every closed path in $Y_{col}$ has exponent sum zero. In particular every closed path in $\partial Y_{col}$ has exponent sum zero. Let $d$ be a 2-cell in $Z_i$ and let $\gamma$ be the attaching path of $d$. Note that $\gamma\subseteq \partial Y_{col}$. Then $f(\gamma)\subseteq K$ is a closed path of exponent sum zero, and therefore there exists a 2-cell $\hat d$ in $\hat K$ with attaching path $f(\gamma)$. Define $f'(d)=\hat d$. This gives the desired extended immersion. 
Note that the essential part of $f'^{-1}(\hat K)$  is 
$$W'=Z_1\cup\ldots Z_p\cup Y_{k+1}\cup\ldots\cup Y_m\cup\ldots\cup Y_n.$$
Since $(\hat L,\hat K)$ satisfies the strong relative coloring test, it follows from Proposition \ref{prop:colim} that $(X',W')$ does so as well. And since the angle structure on $\hat K$ is standard, so is the induced one on $W'$.
Since $Z_i$, $1\le i\le p$, is contractible and all 2-cells are attached along loops of exponent sum zero, there are sink and source vertices in $Z_i$ by Lemma \ref{lem:sink}. Since the angle structure on $Z_i$ is standard, this implies that there is a vertex $v_i\in Z_i$ so that $lk(v_i,Z_i)=lk_0(v_i,Z_i)$. Since $lk(v_i,Z_i)$ is connected by construction this shows that $lk_0(v_i,Z_i)=lk_0(v_i,W')$ is connected. The last equation is true because 
$$W'=Z_1\cup\ldots Z_p\cup Y_{k+1}\cup\ldots Y_m\cup\ldots\cup Y_n$$
is a disjoint union. 
We can now use Lemma \ref{lem:val1} to conclude that $\kappa(v_i,X')<\kappa(v_i,W')$, $1\le i\le p$. 

We know that $Y_i$, $k+1\le i\le m$, is collapsible and $\partial Y_i$ contains a vertex $v_i$ of valency 0 or 1. By Lemma \ref{lem:val1} it follows that $\kappa(v_i,X')< \kappa(v_i,W')$. We have located vertices $v_1,\ldots, v_p,\ldots ,v_m$ for which $\kappa(v_i,X')< \kappa(v_i,W')$.
Thus, by Proposition \ref{prop:heredi} we have 
$$\chi(X')\le \chi(W')-m,$$
and since each $Z_i$, $1\le i\le p$, and each $Y_i$, $k+1\le i \le m$, is collapsible and therefore has Euler characteristic 1, we have 
$$\chi(X')\le \chi(W')-m=\chi(Y_{m+1})+\ldots +\chi(Y_n).$$ Since $Y_i$, $m+1\le i\le n$ is not collapsible and $f\colon Y_i\to K$ is a combinatorial immersion and, furthermore, $K$ has collapsing nonpositive immersion, it follows that $\chi(Y_i)\le 0$. In summary we have
$$\chi(X')\le 0.$$

We have left to show that $$\chi(X)\le \chi(X').$$
$$\chi(X)=\chi(X-Y_{col}^{\circ})+\chi(Y_{col})-
\chi(\partial Y_{col})$$
$$\chi(X')=\chi(X'-Z^{\circ})+\chi(Z)-\chi(\partial Z).$$
Note that $$X'=X_{rest}\cup \mbox{\{ 2-cells of $Z$\}}.$$
As observed above, $Z^{\circ}$ are exactly the open 2-cells in $Z$.
Thus $X'-Z^{\circ}=X_{rest}$. So we have
$$X-Y_{col}^{\circ}=X'-Z^{\circ}=X_{rest}$$ and also
$$\partial Y_{col}=\Delta_1\cup\ldots\cup \Delta_p=\partial Z.$$
Thus 
$$\chi(X)=\chi(X_{rest})+\chi(Y_{col})-(\chi(\Delta_1)+\ldots +\chi(\Delta_p))$$
$$\chi(X')=\chi(X_{rest})+\chi(Z)-(\chi(\Delta_1)+\ldots +\chi(\Delta_p)).$$ We obtain
$$\chi(X)-\chi(X')=\chi(Y_{col})-\chi(Z).$$
Recall that 
$$Y_{col}=Y_1\cup\ldots \cup Y_k$$ 
and the $Y_i$, $1\le i\le k$, are the collapsible components of $Y$ so that $\partial Y_i$ does not contain a vertex of valency 0 or 1,
and
$$Z=Z_1\cup\ldots \cup Z_p \mbox{ (collapsible connected components)}.$$
Therefore $\chi(Y_{col})=k$ and $\chi(Z)=p$. On one hand we have
$$\partial Y_{col}=\partial Y_1\cup\ldots\cup \partial Y_k$$
and on the other 
$$\partial Y_{col}=\Delta_1\cup\ldots\cup \Delta_p \ \mbox{(connected components)}.$$
Thus we have $k\le p$, and it follows that
$$\chi(Y_{col})-\chi(Z)= k-p\le 0.$$ In summary
$$\chi(X)\le \chi(X')\le 0.$$ \qed

\begin{remark}
 In the above proof we have found vertices $v_i\in \Delta_i=\partial Z_i$ so that $\kappa(v_i,X')<\kappa(v_i,W')$. Note that since $\partial Z\subseteq \partial Y$, these vertices are also contained in $Y$.  Here we remark that for these vertices we also have $\kappa(v_i,X)<\kappa(v_i,Y)$. Note that $lk_0(v_i,X')$ is obtained from $lk_0(v_i,X)$ by the following process: The interior of  $lk_0(v_i,Y)$ is removed and replaced with the connected graph $lk_0(v_i,W')$. Since $lk_0(v_i,X')$ contains a connected component $C$ such that $C\cap lk_0(v_i,W')=\emptyset$, $C$ is also a component of $lk_0(v_i,X)$ such that $C\cap lk_0(v_i,Y)=\emptyset$. It follows that $\kappa(v_i,X)<\kappa(v_i,Y)$. This argument shows that every collapsible component of $Y$ contains a vertex $v$ such that $\kappa(v,X)<\kappa(v,Y)$. This provides an alternative proof for the fact that $\chi(X)\le 0$.
 \end{remark}

\section{Labeled oriented trees}

A labeled oriented graph (LOG) $\Gamma = (E, V, s, t, \lambda)$ consists of two sets $E$, $V$ of edges and vertices, and three maps $s, t, \lambda\colon E\to V$ called, respectively source, target and label. $\Gamma$ is said to be a labeled oriented tree (LOT) when the underlying graph is a tree. The associated LOG presentation is defined as
$$P(\Gamma)=\langle V\ |\  s(e)\lambda(e)=\lambda(e)t(e),\ e\in E \rangle.$$
The LOG complex $K(\Gamma)$ is the standard 2-complex defined by the presentation, and the group $G(\Gamma)$ presented by $P(\Gamma )$ is equal to $\pi_1(K(\Gamma))$. 

It is known that LOT-complexes are spines of ribbon 2-knot complements. See Howie \cite{How83}. So the study of LOTs is an extension of classical knot theory. Asphericity, known for classical knots, is unresolved for LOTs. The asphericity question for LOTs is of central importance to Whitehead's asphericity conjecture: A subcomplex of an aspherical 2-complex is aspherical. See Berrick/Hillman \cite{BerrickHillman}, Bogley \cite{Bogley}, and Rosebrock \cite{Ro18}.

A sub-LOG $\Gamma_0=(E_0, V_0)\subseteq\Gamma$ is a subgraph so that $E_0\ne\emptyset$ and $\lambda\colon E_0\to V_0$. A LOG is called {\em boundary reduced} if whenever $v$ is a vertex of valency 1 then $v=\lambda(e)$ for some edge $e$. It is called {\em interior reduced} if for every vertex $v$ no two edges starting or terminating at $v$ carry the same label. It is called {\em compressed }if for every edge $e$ the label $\lambda(e)$ is not equal to $s(e)$ or $t(e)$. Finally, a LOG is {\em reduced} if it is boundary reduced, interior reduced, and compressed. Given a LOG, reductions can be performed to produce a reduced LOG, and, in case the LOG is a LOT, this process does not affect the homotopy type of the LOT complex. A LOG is called {\em injective} if the labeling map $\lambda\colon E\to V$ is injective.

We collect some known facts about injective LOTs (see\\
 Harlander/Rosebrock \cite{HarRose2022})
\begin{itemize}
    \item[F1] (see Theorem 3.2 and Theorem 3.3 of \cite{HarRose2022}) If $\Gamma$ is a reduced injective LOT that does not contain boundary reducible sub-LOTs, then $K(\Gamma)$ admits a zero/one-angle structure that satisfies the coloring test. It follows that $K(\Gamma)$ is DR (and hence aspherical), has collapsing non-positive immersion, and $G(\Gamma)$ is locally indicable.
    \item[F2] (see Theorem 8.4 of \cite{HarRose2022}) Suppose $\Gamma$ is reduced and injective. Assume that $\Gamma_1,\ldots, \Gamma_n$ are disjoint sub-LOTs and that collapsing each into one of its vertices produces a quotient LOT $\bar\Gamma$ without boundary reducible sub-LOTs. Then $K(\Gamma)$ admits a zero/one-angle structure that satisfies a coloring test relative to $K(\Gamma_1),\ldots ;K(\Gamma_n)$. (This relative coloring test agrees with the one defined in this paper only in case $n=1$).
    \item[F3] (see Theorem 8.5 of \cite{HarRose2022}) If $\Gamma$ is reduced and injective, then $K(\Gamma)$ is VA and hence aspherical. It follows that if $\Gamma$ is an injective LOT, reduced or not, then $K(\Gamma)$ is aspherical.
\end{itemize}

We believe that if $\Gamma$ is a reduced injective LOT, then $K(\Gamma)$ has collapsing non-positive immersion. We outline a strategy for proving this. $\Gamma$ contains disjoint maximal sub-LOTs $\Gamma_1,\ldots, \Gamma_n$ so that identifying each into one of its vertices produces a quotient LOT $\bar\Gamma$ without boundary reducible sub-LOTs. Cases that differ from this scenario typically yield to ad hoc considerations. We may assume by induction on the number of vertices that each $K(\Gamma_i)$ has collapsing non-positive immersion. If $n=1$ the methods developed in this paper apply and we obtain a positive result. See Theorem \ref{thm:lotnpi} below. If $n>1$ then we are in the situation of the second point F2 mentioned above. $K(\Gamma)$ does satisfy a relative coloring test, but not the version presented in this paper, and certainly not the strong relative coloring test. Thus the methods developed here fall short in establishing non-positive immersion for all reduced injective LOT complexes $K(\Gamma)$. We do think that our methods can be strengthened to give this result. This is a topic for future work.

\begin{thm}\label{thm:lotnpi}
Suppose $\Gamma$ is a reduced injective LOT that contains a sub-LOT $\Gamma_1$ so that 
\begin{enumerate}
    \item $K(\Gamma_1)$ has collapsing non-positive immersion;
    \item identifying $\Gamma_1$ to one of its vertices produces a reduced quotient LOT $\bar\Gamma$ without boundary reducible sub-LOTs.
\end{enumerate} Then $K(\Gamma)$ has collapsing non-positive immersion.
\end{thm}

\begin{proof}
Let $L=K(\Gamma)$ and $K=K(\Gamma_1)$. Assume that $\Gamma_1$ is collapsed to the vertex $y$. Then $K$ can be folded to the edge $y$ which produces $\bar L=K(\bar\Gamma)$. We will show that the result follows from Theorem \ref{thm:mainapp}.  $\bar L$ satisfies the coloring test by fact F1 stated above. In order to show that $y^+$ and $y^-$ lie in different components of $lk_0(\Bar L)$ we need to take a closer look at that link. It was shown in \cite{HarRose2022}, Theorem 3.3, that $\bar L=K(\bar\Gamma)$ has the following local bi-forest property: If $x_1,\ldots, x_n$ are the edges of $\bar L$, then there exists a choice of $\epsilon_i\in \{ +,-\}$ so that 
$\Lambda_1=\Lambda(x_1^{\epsilon_1},\ldots, x_n^{\epsilon_n})$ and $\Lambda_2=\Lambda(x_1^{-\epsilon_1},\ldots, x_n^{-\epsilon_n})$ are forests. Here $\Lambda(x_1^{\epsilon_1},\ldots, x_n^{\epsilon_n})$ is the subgraph of $lk(\bar L)=\Lambda$ spanned by the vertices $x_1^{\epsilon_1},\ldots, x_n^{\epsilon_n}$. Furthermore, a zero/one-angle structure can be put on $\bar L$ (see \cite{HarRose2022}, Theorem 3.2) so that
$$lk_0(\bar L)=\Lambda_1\cup\Lambda_2.$$ Since $\Lambda_1$ and $\Lambda_2$ are disjoint, it follows that $y^+$ and $y^-$, the $y$ being one of the $x_i$, lie in different components of $lk_0(\bar L)$.
Theorem \ref{thm:mainapp} now gives the result.
\end{proof}

The {\em core} of a LOT $\Gamma$ is the sub-LOT obtained when performing all boundary reductions on $\Gamma$. Note that if the core of $\Gamma_1$ in Theorem \ref{thm:lotnpi} does not contain boundary reducible sub-Lots, then it follows from Lemma \ref{lem:core} and the above fact F1 that $K(\Gamma_1)$ does have collapsing non-positive immersion.

\begin{lemma}\label{lem:core}
Let $\Gamma$ be a LOT and let $\Gamma_c$ be its core. If $K(\Gamma_c)$ has collapsing non-positive immersion, then so does $K(\Gamma)$.
\end{lemma}

\begin{proof}
Let $f\colon X\to K(\Gamma)$ be a combinatorial immersion. Assume that $X$ has no free vertices or edges and is not a point. Then $f$ maps into $K(\Gamma_c)$, because otherwise it would have a free edge. It follows that $\chi(X)\le 0$
\end{proof}

\vspace{3ex}

\noindent Jens Harlander\\
Department of Mathematics\\
Boise State University\\
Boise, ID 83725-1555\\
USA\\
email: jensharlander@boisestate.edu\\

\noindent Stephan Rosebrock\\
P{\"a}dagogische Hochschule Karlsruhe\\
Bismarckstr. 10\\
76133 Karlsruhe\\
Germany\\
email: rosebrock@ph-karlsruhe.de

\end{document}